\documentclass[12pt]{article}

\usepackage{amsmath, amsthm, amssymb, enumerate}
\usepackage{color}
\usepackage{datetime}
\usepackage{fancyhdr}
\usepackage{bbm}
\chardef\coloryes=1 
\chardef\isitdraft=1 

\ifnum\isitdraft=1 
   \def\version{9} 
   \def\eqref#1{({\ref{#1}})}                

\definecolor{labelkey}{gray}{.3}

\definecolor{refkey}{rgb}{.3,0.3,0.3}

\else

  \def\startnewsection#1#2{\section{#1}\label{#2}\setcounter{equation}{0}}   
  \def\nnewpage{} 
\fi

\begin{document}
\def\ques{{\cor \underline{??????}\cob}}
\def\nto#1{{\coC \footnote{\em \coC #1}}}
\def\fractext#1#2{{#1}/{#2}}
\def\fracsm#1#2{{\textstyle{\frac{#1}{#2}}}}   
\def\nnonumber{}


\def\cor{{}}
\def\cog{{}}
\def\cob{{}}
\def\coe{{}}
\def\coA{{}}
\def\coB{{}}
\def\coC{{}}
\def\coD{{}}
\def\coE{{}}
\def\coF{{}}

\ifnum\coloryes=1

  \definecolor{coloraaaa}{rgb}{0.1,0.2,0.8}
  \definecolor{colorbbbb}{rgb}{0.1,0.7,0.1}
  \definecolor{colorcccc}{rgb}{0.8,0.3,0.9}
  \definecolor{colordddd}{rgb}{0.0,.5,0.0}
  \definecolor{coloreeee}{rgb}{0.8,0.3,0.9}
  \definecolor{colorffff}{rgb}{0.8,0.3,0.9}
  \definecolor{colorgggg}{rgb}{0.5,0.0,0.4}

 \def\cog{\color{colordddd}}
 \def\cob{\color{black}}
 \def\cor{\color{red}}
 \def\coe{\color{colorgggg}}

 \def\coA{\color{coloraaaa}}
 \def\coB{\color{colorbbbb}}
 \def\coC{\color{colorcccc}}
 \def\coD{\color{colordddd}}
 \def\coE{\color{coloreeee}}
 \def\coF{\color{colorffff}}
 \def\coG{\color{colorgggg}}

\fi
\ifnum\isitdraft=1
   \chardef\coloryes=1 
   \baselineskip=17pt
\pagestyle{myheadings}
\reversemarginpar

\def\const{\mathop{\rm const}\nolimits}  
\def\diam{\mathop{\rm diam}\nolimits}    

 \def\llabel#1{\label{#1}{\ \mbox{\rm\color{red} {#1}\color{black}}}}

\def\rref#1{{\ref{#1}{\rm \tiny \fbox{\tiny #1}}}}
\def\theequation{\fbox{\bf \thesection.\arabic{equation}}}
\def\ccite#1{{\cite{#1}{\rm \tiny ({#1})}}}

\def\startnewsection#1#2{\newpage\cog \section{#1}\cob\label{#2}

\setcounter{equation}{0}
\pagestyle{fancy}

\lhead{\cob Section~\ref{#2}, #1 }
\cfoot{}
\rfoot{\thepage}
\lfoot{\cob{\today,~\currenttime}~(c75-iklt2, Version~\fbox{\version})}}
\chead{}
\rhead{\thepage}
\def\nnewpage{\newpage}

\newcounter{startcurrpage}
\newcounter{currpage}

\def\llll#1{{\rm\tiny\fbox{#1}}}
   \def\blackdot{{\color{red}{\hskip-.0truecm\rule[-1mm]{4mm}{4mm}\hskip.2truecm}}\hskip-.3truecm}
   \def\bdot{{\coC {\hskip-.0truecm\rule[-1mm]{4mm}{4mm}\hskip.2truecm}}\hskip-.3truecm}
   \def\purpledot{{\coA{\rule[0mm]{4mm}{4mm}}\cob}}
   \def\pdot{\purpledot}
\else
   \baselineskip=15pt
   \def\blackdot{{\rule[-3mm]{8mm}{8mm}}}
   \def\purpledot{{\rule[-3mm]{8mm}{8mm}}}
   \def\pdot{}
\fi

\def\tdot{\fbox{\fbox{\bf\tiny I'm here; \today \ \currenttime}}}
\def\nts#1{{\hbox{\bf ~#1~}}} 
\def\nts#1{{\cor\hbox{\bf ~#1~}}} 
\def\ntsf#1{\footnote{\hbox{\bf ~#1~}}} 
\def\ntsf#1{\footnote{\cor\hbox{\bf ~#1~}}} 
\def\bigline#1{~\\\hskip2truecm~~~~{#1}{#1}{#1}{#1}{#1}{#1}{#1}{#1}{#1}{#1}{#1}{#1}{#1}{#1}{#1}{#1}{#1}{#1}{#1}{#1}{#1}\\}
\def\biglineb{\bigline{$\downarrow\,$ $\downarrow\,$}}
\def\biglinem{\bigline{---}}
\def\biglinee{\bigline{$\uparrow\,$ $\uparrow\,$}}

\newtheorem{Theorem}{Theorem}[section]
\newtheorem{Corollary}[Theorem]{Corollary}
\newtheorem{Proposition}[Theorem]{Proposition}
\newtheorem{Lemma}[Theorem]{Lemma}
\newtheorem{Remark}[Theorem]{Remark}
\newtheorem{Example}[Theorem]{Example}
\newtheorem{definition}{Definition}[section]
\def\theequation{\thesection.\arabic{equation}}
\def\endproof{\hfill$\Box$\\}
\def\square{\hfill$\Box$\\}
\def\comma{ {\rm ,\qquad{}} }            
\def\commaone{ {\rm ,\qquad{}} }         
\def\dist{\mathop{\rm dist}\nolimits}    
\def\sgn{\mathop{\rm sgn\,}\nolimits}    
\def\Tr{\mathop{\rm Tr}\nolimits}    
\def\div{\mathop{\rm div}\nolimits}    
\def\supp{\mathop{\rm supp}\nolimits}    
\def\divtwo{\mathop{{\rm div}_2\,}\nolimits}    
\def\re{\mathop{\rm {\mathbb R}e}\nolimits}    
\def\div{\mathop{\rm{Lip}}\nolimits}   
\def\indeq{\qquad{}}                     
\def\period{.}                           
\def\semicolon{\,;}                      

\title{Continuity of Cost Functional and Optimal Feedback Controls for the Stochastic Navier Stokes Equation in 2D}
\author{Kerem U\u{g}urlu}
\maketitle

\date{}

\begin{center}
\end{center}

\medskip

\indent Department of Mathematics, University of Southern California, Los Angeles, CA 90089\\
\indent e-mail:kugurlu@usc.edu

\begin{abstract}
We show the continuity of a specific cost functional  $J(\phi) = \mathbb{E} \sup_{ t \in  [0,T]}(\varphi(\mathcal{L}[t,u_\phi(t), \phi(t)]))$  of the SNSE in 2D on an open bounded nonperiodic domain $\mathcal{O}$ with respect to a special set of feedback controls $\{\phi_n\}_{n \geq 0}$, where $\varphi(x) =\log(1 + x)^{1-\epsilon}$ with $0 < \epsilon <1$. 
\end{abstract}


\noindent
\section{Introduction}

We consider the existence of optimal controls $\phi$ related to the stochastic Navier-Stokes equation (SNSE) describing the flow of a viscous incompressible fluid in a smooth bounded domain $\mathcal{O} \subset \mathbb{R}^2$ with a multiplicative white noise 
\begin{align}
\label{main_equation}
du_\phi + (u_\phi \cdot \nabla u_\phi  + \nabla p - \nu \Delta u_\phi - \phi) dt &= g(u_\phi)d\mathcal{W},
\nonumber \\
\nabla \cdot u_\phi &= 0,
\nonumber \\
u(0) &= u_0
\end{align}
and Dirichlet boundary condition $u_\phi = 0$ on $[0,\infty) \times \partial \mathcal{O}$, where $u_\phi = (u_{\phi_1}, u_{\phi_2})$ represent the velocity field, $\nu$ stands for the coefficient of kinematic viscosity, $\phi$ is the deterministic force and $p$ represents the pressure. $g(u_\phi)d\mathcal{W} = \sum_{k=1}^\infty g_k(u_\phi)dW_k$ is the cylindrical Brownian motion with independent one dimensional Brownian motions $W_k$ and Lipschitz coefficients $g_k(u_\phi)$ \cite{BKL00,BP00,BT73,CG94,C89,CP97,DD03,FG95,FR02,GV14,M02, K06,MR05,MS02,O06,S03}. Compared to their deterministic counterparts, the stochastic PDE's lead us to consider new questions and additional technical difficulties such as the existence and uniqueness of invariant measures or lack of compactness caused by the stochastic term driven by Brownian motion. There are two notions for the solutions of the SNSE. The first notion is called a \textit{martingale} solution, where the stochastic basis is not given in advance but constructed as a part of the solution \cite{CG94, C89, FG95, MR04, V76}. The second notion, which we consider here, is the so called \textit{pathwise} solution or \textit{strong} solution, i.e. a complete probability space and the Brownian motion are given a priori \cite{B00,GV14,GZ09,KV}.

The existence of optimal control of stochastic evolution equations has been studied by \cite{PI85,G93,G87,GS94,T89,T90} among others  by adding linearity or semilinearity assumptions as well as putting boundedness restriction for nonlinearities. On the other hand, more specifically, the literature about the optimal control of SNSE is not very rich, since these assumptions do not apply to the SNSE. The nonlinearity of the SNSE causes the problem to be of non-convex type. We refer the reader the book of Ekeland and Temam about non-convex optimization for further investigations \cite{ET99}. Related works about the control of SNSE can be mentioned as follows. Choi et al. investigated the optimal control problem in \cite{C93} for the stochastic Burgers equation  (one dimensional Navier-Stokes equation) with additive noise. The paper \cite{PD99} studies the control of turbulence for the stochastic Burgers equation. Another work in this direction is by Sritharan \cite{S00}, where the existence of optimal controls is established using techniques for the martingale problem formulation of Stroock and Varadhan \cite{SV79} in the context of stochastic Navier-Stokes equation. In \cite{B00}, it is shown that there exist feedback controls for the SNSE of \eqref{main_equation}, which are controlled by different external forces $\phi$ for a specific cost functional satisfying some regularity conditions. In our paper, we follow the framework that is studied in \cite{B00}. We show using the recent bounds and approximations for the SNSE in 2D in \cite{KUZ} that the cost functional $J(\phi)$ is continuous with respect to a specific set $\phi$. Contrary to \cite{B00} though, we control the \textit{supremum} of SNSE up to a terminal deterministic time $T$, which is only natural to introduce when we want to control the extreme events on the whole path rather than integrating the path.

The rest of the paper is as follows: in Section 2, we give the functional setting, which describes the assumptions on the problem, deterministic and stochastic framework as well as the notion of the solution considered throughout the paper. In Section 3, we state and prove our main result for the control problem. Section 4 gives the crucial technical results used in our result.
\section{Functional Setting}
First, we recall the deterministic and probabilistic framework used throughout the paper. 
\subsection{Deterministic Framework}
Let $\mathcal{O}$ be a bounded open connected subset of $\mathbb{R}^2$ and $\partial\mathcal{O}$ be smooth. We take \\ $\mathcal{V} = \{u \in C_0^{\infty}(\mathcal{O})^2: \nabla \cdot u = 0 \}$ and denote by $H$ the closure of $\mathcal{V}$ in $L^2(\mathcal{O})$ and $V$ the closure of $\mathcal{V}$ in $H^1(\mathcal{O})$, respectively. Hence, the spaces $H$ and $V$ are identified by
\begin{align}
\label{v_space}
H &= \{ u \in L^2(\mathcal{O})^2: \nabla \cdot u = 0, u \cdot n |_{\partial \mathcal{O}} = 0 \}, \\
V &= \{ u \in H_0^1(\mathcal{O})^2: \nabla \cdot u = 0 \}. 
\end{align}
Here $n$ is the outer pointing normal to $\partial \mathcal{O}$. On $H$ we take the $L^2(\mathcal{O})$ inner product and norm  as 
\begin{equation*}
\langle u,v \rangle := \int_{\mathcal{O}} u \cdot v d\mathcal{O}, \qquad \lVert u \rVert_H := \sqrt{\langle u,u \rangle}.
\end{equation*}
We denote the inner product on $H$ by $\langle \cdot, \cdot \rangle$ and the norm by $\lVert \cdot \rVert_{H}$. The Leray-Hopf projector, $P_H$ is defined as the orthogonal projection of $L^2(\mathcal{O})^2$ onto $H$. Moreover, on $V$, we use the $H^1$ norm and inner products
\begin{equation}
\langle\langle u,v \rangle\rangle := \int_{\mathcal{O}} \nabla u \cdot \nabla v d\mathcal{O}, \qquad \lVert u \rVert_V := \sqrt{\langle\langle u,u \rangle\rangle }, \qquad u,v \in V.
\end{equation}
We note here that due to Dirichlet boundary condition in Equation \ref{main_equation}, the Poincare inequality 
\begin{equation}
\lVert u \rVert_H \leq C \lVert u \rVert_V, \qquad \forall u \in V
\end{equation}
holds, justifying $\lVert \cdot \rVert_V$ as a norm. We take $V'$ to be the dual of $V$, relative to $H$ with the pairing notated by $\langle \cdot, \cdot \rangle$.
We  next define the Stokes operator $A$. $A$ is understood as a bounded linear map from $V$ to $V'$ via: 
\begin{equation}
\langle Au, v \rangle = \langle \langle u,v \rangle \rangle \qquad u,v \in V.
\end{equation}
 $A$ can be extended to an unbounded operator from $H$ to $H$ according to $Au = -P_H\Delta u$ with the domain $\mathcal{D}(A)= V \cap H^2(\mathcal{O})$. The dual of $V = \mathcal{D}(A^{1/2})$ with respect to $H$ is denoted by $V'= \mathcal{D}(A^{-1/2})$.
By the theory of symmetric, compact operators for $A^{-1}$, we have the existence of an orthonormal basis $\{e_k\}$ for $H$ of eigenfunctions of $A$. We recall here that the corresponding eigenvalues $\{ \lambda_k \}$ form an increasing, unbounded sequence of 
\begin{equation*}
0 < \lambda_1 \leq \lambda_2 \leq ... \leq \lambda_n \leq ...
\end{equation*}
We also define the nonlinear term as a bilinear mapping $ V \times V$ to $V'$ via 
\begin{equation*}
B(u,v)= \mathcal{P}_H(u \cdot \nabla v )
\end{equation*} 
We note here that the cancellation property $\langle B(u,v),v \rangle = 0$ holds for $u,v \in V$. Moreover, we denote by $\mathcal{L}(H)$ and $\mathcal{L}(V)$ as the space of all linear and continuous operators from the Banach space $H$ and $V$ to themselves respectively. 
\subsection{Stochastic Framework}
In this section, we recall the necessary background material in stochastic analysis in infinite dimensions needed in the paper (see \cite{DZ92, DGT11, F08, PR07}). We fix a stochastic basis $\mathcal{S}= ( \Omega, \mathcal{F}, \mathbb{P}, \{ \mathcal{F}_t\}, \mathcal{W} ) $, which consists of a complete probability space ($\Omega,\mathbb{P}$), equipped with a complete right-continuous filtration $\mathcal{F}_t$, and a cylindrical Brownian motion $\mathcal{W}$, defined on a separable Hilbert space $U$ adapted to this filtration. 

Given a separable Hilbert space $X$, we denote by $L_2(U,X)$ the space of Hilbert-Schmidt operators from $U$ to $X$, equipped with the norm $\lVert G \rVert_{L_2(U,X)}= (\sum_k \lVert G \rVert^2_X)^{1/2}$ \cite{DZ92}. For an $X$-valued predictable process 
$G \in L^2( \Omega;L^2_{\text{loc}}([0,\infty]);L_2(U,X) )$, we define the Ito stochastic integral

\begin{equation}
\label{def1}
\int_0^t G d\mathcal{W} := \sum_k \int_0^t G_k dW_k
\end{equation}
which lies in the space $\mathcal{O}_X$ of $X$-valued square integrable martingales. We also recall the Burkholder-Davis-Gundy inequality: For any $p \geq 1$ we have 
\begin{equation}
\label{def2}
\mathbb{E} \left ( \sup_{t \in [0,T]} \rVert \int_0^t Gd\mathcal{W} \lVert^p_X \right ) \leq C\mathbb{E} \left ( \int_0^T  \rVert  G \lVert^2_{L_2(U,X)} \right )^{p/2} 
\end{equation}
for some $C=C(p) > 0$.
\subsection{Conditions on the Noise}
Given a pair of Banach spaces $X$ and $Y$, we denote by $\rm{Lip}_u(X,Y)$ the collection of continuous functions $h \colon [0,\infty) \times X \rightarrow Y$ which are sublinear 
\begin{equation}
\label{def3}
\lVert h(t,x) \rVert_Y \leq K_Y(1 + \rVert x \lVert_X)\comma t \geq 0, x \in X  
\end{equation}
and Lipschitz 
\begin{equation}
\label{def4}
\lVert h(t,x) - h(t,y) \rVert_Y \leq K_Y \lVert x - y \rVert_X \comma t \geq 0\commaone x,y \in X
\end{equation}
for some constant $K_Y > 0$ independent of $t$.
The noise term $g(u)d\mathcal{W}$ is defined by
\begin{equation}
\label{def5}
g = \{g_k\}_{k \geq 1}\colon [0,\infty) \times H \rightarrow L_2(U,H).
\end{equation}
Namely
\begin{equation}
\lVert g(t,x) \rVert_{L_2(U, H)} \leq K_j(1 + \lVert x \rVert_H)
\end{equation}
and
\begin{align} 
\lVert g(t,x) - g(t,y) \rVert_{L_2(U,H)} &\leq K_1\lVert x  - y \rVert_H
\nonumber \\
\lVert g(t,x) - g(t,y) \rVert_{L_2(U,V)} &\leq K_2 \lVert x  - y \rVert_V
\end{align}
with 
\begin{equation}
g \in \rm{Lip}_u(H, L^2(U,H)) \cap \rm{Lip}_u(V,L_2(U,V)).
\end{equation}
Given $u \in L^2(\Omega; L^2(0,T;H))$ and $g$ as above the stochastic integral $\int_0^t g(u)d\mathcal{W}$ is a well-defined $H$-valued Ito stochastic integral that is predictable and is such that 
\[
\left \langle \int_0^t g(u)d\mathcal{W}, v \right \rangle = \sum_k \int_0^t \langle g_k(u),v \rangle d\mathcal{W}_k
\]
holds for any $v \in H$.

\subsection{Notion of Solution}
We consider \emph{strong pathwise solutions}, which are solutions with values in $V$ and strong in the \emph{probabilistic} sense, i.e., the driving noise and the filtration are given in advance.
\begin{definition} \label{global_def}
We fix a stochastic basis $\mathcal{S}$ and g is as above. We further assume that the initial data $u_0 \in L^4(\Omega;H) \cap L^2(\Omega;V)$ is $\mathcal{F}_0$ measurable. Moreover, we take $ \phi \subset \mathcal{U}$ of bounded feedback controls. Namely, for each time $t \in T$, we have $\phi(t,.)$ is a continuous linear functional from $H$ to $H$. Moreover, we assume that
\begin{align}
\label{f_assum}
\sup_{t \in [0,T]} \lVert \phi(t,\omega) \rVert_V &\leq K, a.s.,
\nonumber \\ 
\lVert \phi(t_1,x_1) - \phi(t_2,x_2) \rVert_V^2 &\leq C_1 \vert t_1 -t_2 \rvert^2 + C_2 \lVert x_1 - x_2 \rVert^2_V,
\nonumber \\
\lVert \phi(t_1,x_1) - \phi(t_2,x_2) \rVert_H^2 &\leq C_1 \vert t_1 -t_2 \rvert^2 + C_2 \lVert x_1 - x_2 \rVert^2_H
\end{align}
where $C_1$ and $C_2$ are uniform for the family of controls $\phi \in \mathcal{U}$. Then, we say that the pair $(u_\phi,\tau)$ is called a pathwise strong solution of the system if $\tau$ is a stricly positive stopping time, $u_\phi(\cdot \wedge \tau )$ is a predictable process in H such that 
\begin{equation}
\label{init_ass}
u_\phi(\cdot \wedge \tau ) \in L^2(\Omega;C([0,\infty);V))
\end{equation}
with
\begin{equation}
u_\phi\mathbf{1}_{t\leq \tau} \in L^2(\Omega;C([0,\infty);\mathcal{D}(A)))
\end{equation}
and if
\begin{equation}
\label{main_eqn}
\langle u_\phi(t \wedge \tau), v \rangle + \int_0^{t \wedge \tau} \langle \nu Au_\phi + B(u_\phi,u_\phi) - \phi, v \rangle dt = \langle u_0,v\rangle + \sum_k\int_0^{t \wedge \tau} \langle g_k(u_\phi),v \rangle dW_k
\end{equation}
holds for every $v \in H$. Moreover, $(u,\xi)$ is called a maximal pathwise strong solution if $\xi$ is a strictly positive stopping time and there exists $\tau_n \rightarrow \xi$ increasing such that (u, $\tau_n$) is a local strong solution and 
\begin{equation}
\label{stop_time}
\sup_{t \in [0,\tau_n]} \lVert u_\phi \rVert^2_V + \nu\int_0^{\tau_n} \lVert Au_\phi \rVert^2_H dt \geq n
\end{equation}
on the set $\{ \xi < \infty \}$. Such a solution is called global if $\mathbb{P}(\xi < \infty ) = 0$.
\end{definition}
\section{Main Result}\label{Rationality}
Our purpose is to control the solution $u_{\phi}$ of the SNSE 
\begin{align} 
\label{main_eqn}
&\langle u_\phi(t),v \rangle + \int_0^t \langle Au_\phi(s),v \rangle ds 
\nonumber\\&\indeq
= \langle u_0,v \rangle + \int_0^t\langle B(u_\phi(s), u_\phi(s)), v \rangle ds + \int_0^t \langle \phi(s,u_\phi),v \rangle ds 
\nonumber\\&\indeq \indeq 
+ \int_0^t \langle g(s,u_\phi(s)),v \rangle dW_s,
\end{align}
for all $v \in H$, $t \in [0,T]$, a.e. $\omega \in \Omega$, where $\phi \in \mathcal{U}$ and the assumptions on initial data are as in Definition 2.1. We consider the following cost functional 
\begin{equation} 
\label{cfdef}
J(\phi) = \mathbb{E} \sup_{ t \in  [0,T]}( \varphi(\mathcal{L}[t,u_\phi(t), \phi(t)])), 
\end{equation}
with the objective of minimization of the cost functional
\begin{equation}
\min_{\phi \in \mathcal{U}} J(\phi),
\end{equation}
where $\mathcal{L}: [0,T] \times V \times H \rightarrow \mathbb{R}_+$ is Lipschitz i.e. 
\begin{equation}
\lvert \mathcal{L}(t,x_1,y_1)  -  \mathcal{L}(t,x_2,y_2) \rvert^2 \leq C( \lVert x_1 - x_2 \rVert_V^2 + \lVert y_1 - y_2 \rVert_H^2 ),
\end{equation} 
with 
\begin{equation}  
\varphi(x) = \log( 1 + x)^{1-\epsilon},
\end{equation}
and $0 < \epsilon <1$.
\begin{Remark} Since we require $\mathcal{L}$ to be only uniformly Lipschitz, using the concave function $\varphi(x) = \log(1+x)^{1-\epsilon}$ does not imply that, we should check only the end points for the functional $\mathcal{L}$. We also note that the results in the paper still hold for the controls with the support in open subsets of the domain. In their paper F. Abergel and R. Temam \cite{AT90} investigate the deterministic Navier-Stokes equation by controlling the turbulence inside the flow. They give a cost functional regarding the vorticity in the fluid. For our problem, the correspondent functional $\mathcal{L}$ would be 
\begin{equation}  
\mathcal{L}(t,u_\phi(t),\phi(t)) := \lVert \nabla \times u_\phi(t) \rVert_H
\end{equation}
\end{Remark}
We state now our main result. 
\begin{Theorem} \label{thm1} Let $\{\phi_n\}_{n\geq 0}$ be a sequence in $\mathcal{U}$ as defined in \ref{f_assum} and $J$ is as defined in \ref{cfdef}.  Suppose 
\begin{equation}
\sup_{t \in [0,T]} \lVert \phi - \phi_n \rVert_{\mathcal{L}(V)} \rightarrow 0 \textrm{ a.s. }
\end{equation}
Then, we have
\begin{equation} 
J(\phi^n) \rightarrow J(\phi)
\end{equation}
as $n \rightarrow \infty$.
\end{Theorem}
To prove Theorem \ref{thm1}, we use the following result from \cite{B00}.
\begin{Theorem}\cite{B00} \label{breckner} Let $\{\phi_n\}_{n\geq 0}$ be a sequence in $\mathcal{U}$ as defined in  \ref{f_assum}  and suppose 
\begin{equation}
\mathbb{E}[\int_0^T \lVert \phi - \phi_n \rVert_{\mathcal{L}(H)}^2] \rightarrow 0,
\end{equation}
then it holds that 
\begin{equation}
\mathbb{E}[\int_0^T \lVert u_\phi - u_{\phi_n} \rVert_V^2] \rightarrow 0,
\end{equation}
as $n \rightarrow \infty$.
\end{Theorem}
Next, we need the two technical lemmas. 
\begin{Lemma} Let $\varphi(x)$ be a concave increasing function with $\varphi(0) = 0$. Then, we have
\begin{equation}
|\varphi(x_1) - \varphi(x_2) | \leq \varphi(|x_1 - x_2|)
\end{equation}
\end{Lemma}
\begin{proof} 
Since $\varphi(x)$ is a concave increasing function with $\varphi(0) = 0$, $\varphi(x)$ is subadditive. Hence, we have 
\begin{align}
\varphi(x_1) &\leq \varphi( |x_1 - x_2| + x_2 )
\nonumber \\
\varphi( |x_1 - x_2| + x_2 ) &\leq \varphi( |x_1 - x_2|) + \varphi(x_2)
\nonumber \\
\varphi( x_1 ) - \varphi(x_2) &\leq \varphi( |x_1 - x_2|)
\end{align}
By interchanging $x_1$ and $x_2$ we conclude the proof.  
\end{proof}
\begin{Lemma} \cite{GZ09} \label{gronwall} Fix $T > 0$. Assume that 
\begin{equation}
X,Y,Z,R : [0,T) \times \Omega \rightarrow \mathbb{R} 
\end{equation}
are real-valued, non-negative stochastic processes. Let $\tau < T$ be a stopping time so that 
\begin{equation}
\mathbb{E}\int_0^\tau (RX + Z) ds < \infty.
\end{equation}
Assume, moreover that for some fixed constant $\kappa$ we have
\begin{equation}
\int_0^\tau R ds < \kappa, \textrm{ a.s. }
\end{equation}
Suppose that for all stopping times $0 \leq \tau_a \leq \tau_b \leq \tau$
\begin{equation}
\label{eqn1}
\mathbb{E}\bigg( \sup_{t \in [\tau_a,\tau_b]} X + \int_0^\tau Y ds \bigg) \leq C_0 \mathbb{E}\bigg( X(\tau_a) + \int_{\tau_a}^{\tau_b} (RX + Z)ds \bigg),
\end{equation}
where $C_0$ is a constant independent of the choice of $\tau_a,\tau_b$. Then we have 
\begin{equation}
\mathbb{E}\bigg( \sup_{t \in [0,\tau]} X + \int_0^\tau Y ds \bigg) \leq C \mathbb{E}\bigg( X(0) + \int_0^\tau Z ds \bigg), 
\end{equation}
where $C$ depends on $C_0, T$ and $\kappa$.
\end{Lemma}
\begin{proof} Choose a finite sequence of stopping times 
\begin{equation}
 0 = \tau_0 < \tau_1 < ... < \tau_N < \tau_{N+1} = \tau
\end{equation}
so that 
\begin{equation}
\label{eqn2}
\int_{\tau_{k-1}}^{\tau_k} R ds < \frac{1}{2C_0}\textrm{ a.s. }
\end{equation}
For each pair $\tau_{k-1}, \tau_k$ take $\tau_a = \tau_{k-1}$ and $\tau_b = \tau_k$ in (\ref{eqn1}). Using (\ref{eqn2}), we have 
\begin{equation}
\mathbb{E}\bigg( \sup_{t \in [\tau_{k-1},]} X + \int_{\tau_{k-1}}^{\tau_k} Yds \bigg) \leq C \mathbb{E}X(\tau_{k-1}) + C\mathbb{E}\int_{\tau_{k-1}}^{\tau_k}Zds.
\end{equation}
By induction we have 
\begin{equation}
\mathbb{E}\bigg( \sup_{t \in [0,\tau_j]}X + \int_0^{\tau_j}Yds \bigg) \leq C \mathbb{E}X(0) + C \mathbb{E}\int_0^{\tau_j}Zds
\end{equation}
then we have
\begin{align}
\mathbb{E}\bigg( \sup_{t \in [0,\tau_{j+1}]}X + \int_0^{\tau_j}Yds \bigg) &\leq C \mathbb{E}X(0) + C \mathbb{E}\int_0^{\tau_j}Zds + C\mathbb{E}\bigg( \sup_{t \in [\tau_j,\tau_{j+1}]}X + \int_{\tau_j}^{\tau_{j+1}}Yds \bigg)
\nonumber \\&
\leq C \mathbb{E}X(0) + C\mathbb{E}\int_0^{\tau_{j+1}}Zds + C\mathbb{E}X(\tau_j)
\nonumber \\&
\leq C\mathbb{E}X(0) + C\mathbb{E}\int_0^{\tau_{j+1}}Zds
\end{align}
Hence, we conclude the proof. 
\end{proof}
We continue with the following theorem. 
\begin{Theorem} \label{main_prop}
Let $\tilde{M} > 0$ and $M > 1$. Moreover, let $\{u_{\phi_n}\}_{n \geq 1}$ be the sequence of solutions of (3.19). Suppose we have
\begin{equation}
\lVert u_0 \rVert_V \leq \tilde{M}, \textrm{ a.s. }
\end{equation} 
Moreover, assume 
\begin{equation}
\sup_{t \in [0,T]} \lVert \phi - \phi_n \rVert_{\mathcal{L}(V)} \rightarrow 0 \textrm{ a.s. }
\end{equation}
as $m,n \rightarrow \infty$.
Denote 
\begin{equation}
\mathcal{T}_n^{M,T} = \{ \tau \leq T : ( \sup_{t \in [0,\tau]} \lVert u_{\phi_n} \rVert_V^2 + \int_0^\tau \lvert Au_{\phi_n} \rvert^2 dt )^{1/2} \leq M + \tilde{M} \}. 
\end{equation}
Let 
\begin{equation}
\mathcal{T}_{m,n}^{M,T} := \mathcal{T}_m^{M,T} \cap  \mathcal{T}_n^{M,T}.
\end{equation}
Then 
\begin{enumerate}
\item For any $T > 0$, we have 
\begin{equation}
\label{cond1}
\lim_{n \rightarrow \infty } \sup_{m \geq n} \sup_{\tau \in \mathcal{T}_{m,n}^{M,T}} \mathbb{E}[\sup_{t \in [0,\tau]} \lVert u_{\phi_m} - u_{\phi_n} \rVert_V^2 + \int_0^{\tau} \lVert A(u_{\phi_m} - u_{\phi_n}) \rVert_H^2dt  ] = 0.
\end{equation}
\item 
\begin{equation}
\label{cond2}
\lim_{S \rightarrow 0} \sup_n \sup_{\tau \in \mathcal{T}_n^{M,T}} \mathbb{P} \bigg( \sup_{t \in [0,\tau \wedge S]} \lVert u_{\phi_n} \rVert_V^2 
+ \int_0^{\tau \wedge S} \lVert  A(u_{\phi_n} \rVert_H^2 dt > \tilde{M}^2 + (M-1)^2 \bigg) = 0.
\end{equation}
\end{enumerate}
\end{Theorem}
\begin{proof}
1. We have 
\begin{align} 
d(u_{\phi_n} - u_{\phi_m}) + A(u_{\phi_n} - u_{\phi_m})dt &= (B(u_{\phi_n}) - B(u_{\phi_m}))dt + \sum_{k=1}^\infty [g_k(u_{\phi_n}) - g_k(u_{\phi_m})]dW_k 
\nonumber \\&\indeq
+ \phi_n(t,u_{\phi_n}(t)) - \phi_m(t,u_{\phi_m}(t))dt
\end{align}
Hence by Ito-lemma, we have
\begin{align}
&d\lVert u_{\phi_n} - u_{\phi_m} \rVert_V^2 + 2 \lVert A(u_{\phi_n} - u_{\phi_m}) \rVert_H^2 dt 
\nonumber \\&\indeq
= 2 \langle B(u_{\phi_n}) - B(u_{\phi_m}), A(u_{\phi_n} - u_{\phi_m}) \rangle dt 
\nonumber \\&\indeq \indeq
+ \sum_{k=1}^\infty \lVert g_k(u_{\phi_n}) - g_k(u_{\phi_m}) \rVert_V^2dt 
\nonumber \\&\indeq \indeq
+ 2 \sum_{k=1}^\infty \langle g_k(u_{\phi_n}) - g_k(u_{\phi_m}) \rangle dW_k 
\nonumber \\&\indeq \indeq
+ 2 \langle \phi_n(t,u_{\phi_n}) - \phi_m(t,u_{\phi_m}), A( u_{\phi_n} - u_{\phi_m} ) \rangle dt
\end{align}
By taking supremum up to $\tau$, integrating and taking expectation, we get 
\begin{align}
\label{main_tau}
&\mathbb{E}[\sup_{[0,\tau]} \lVert u_{\phi_n} - u_{\phi_m} \rVert_V^2 + 2\int_0^\tau \lVert A(u_{\phi_n} - u_{\phi_m}) \rVert_H^2 dt ] 
\nonumber \\&\indeq
\leq 2 \mathbb{E}\int_0^\tau | \langle B(u_{\phi_m}) - B(u_{\phi_n}), A(u_{\phi_n} - u_{\phi_m}) \rangle |dt 
\nonumber \\&\indeq \indeq
+ \mathbb{E} \int_0^\tau \sum_{k=1}^\infty \lVert g_k(u_{u_{\phi_m} }) - g_k(u_{\phi_n} ) \rVert_V^2dt 
\nonumber \\&\indeq \indeq
+ \mathbb{E} [ \sup_{r \in [0,\tau]} | 2 \sum_{k=1}^\infty \int_0^r \langle g_k(u_{\phi_m}) - g_k(u_{\phi_n}), A(u_{\phi_n} - u_{\phi_m}) \rangle dW_k | ]
\nonumber \\&\indeq \indeq
+ \mathbb{E}[ 2 \int_0^\tau | \langle \phi_n (t, u_{\phi_n}) - \phi_m (t, u_{\phi_m}), A(u_{\phi_n} - u_{\phi_m})  \rangle | dt  ]
\end{align}
We treat each term above seperately. First 
\begin{align}
&\mathbb{E}\int_0^{\tau} \sum_{k=0}^\infty \lVert g_k(u_{\phi_n}) - g_k(u_{\phi_m}) \rVert_V^2 dt 
\nonumber \\& \indeq
\leq \mathbb{E}\int_0^{\tau} \lVert u_{\phi_n} - u_{\phi_m} \rVert_V^2
\end{align}
By Poincare lemma, assumption 3.43 implies that
\begin{equation}
    \int_0^T \lVert u_{\phi_n} - u_{\phi_m}  \rVert_H^2 \rightarrow 0,
\end{equation}
as $m,n \rightarrow \infty$. Then this implies by Theorem \ref{breckner} 
\begin{equation}
\mathbb{E}\int_0^T \lVert u_{\phi_n} - u_{\phi_m} \rVert_V^2 \rightarrow 0,
\end{equation}
as $n,m \rightarrow \infty$.
Next we have 
\begin{align}
&\mathbb{E}[ 2 \int_0^\tau | \langle \phi_n (t, u_{\phi_n}) - \phi_m (t, u_{\phi_m}), A(u_{\phi_n} - u_{\phi_m})  \rangle | dt ]
\nonumber \\& \indeq 
\leq C \mathbb{E}\int_0^T \lVert u_{\phi_n} - u_{\phi_m} \rVert_V^2, 
\end{align}
which goes to 0 as $n,m \rightarrow \infty$ by Theorem \ref{breckner}.

Next, we treat the nonlinear term by seperating into 2 parts as follows.  
\begin{align*}
&| \langle B(u_{\phi_m}- u_{\phi_n}), A(u_{\phi_n} - u_{\phi_m}) \rangle | \leq 
|\langle B(u_{\phi_m}- u_{\phi_n},u_{\phi_m}), A(u_{\phi_n} - u_{\phi_m}) \rangle |
\nonumber \\ \indeq \indeq &+ 
| \langle B(u_{\phi_n}, u_{\phi_n} - u_{\phi_m}), A(u_{\phi_n} - u_{\phi_m}) \rangle |
\end{align*}
For the first term above, we have
\begin{align}
& | \langle B(u_{\phi_m} - u_{\phi_n}, u_{\phi_m}), A(u_{\phi_n} - u_{\phi_m}) \rangle |
\nonumber \\& \indeq 
\leq \lVert u_{\phi_n} - u_{\phi_m} \rVert_V \lVert u_{\phi_m} \rVert_V^{1/2}\lVert Au_{\phi_m} \rVert_H^{1/2}\lVert A(u_{\phi_n} - u_{\phi_m}) \rVert_H
\nonumber \\& \indeq 
\leq \frac{1}{6}\nu \lVert A(u_{\phi_n} - u_{\phi_m})  \rVert_H^2 + C \lVert u_{\phi_n} - u_{\phi_m} \rVert_V^2 \lVert Au_{\phi_m} \rVert_H,
\end{align}
To estimate the term  C $\lVert u_{\phi_n} - u_{\phi_m} \rVert_V^2 \lVert Au_{\phi_m} \rVert_H$ in the second line of inequality, we apply Theorem \ref{gronwall} with $R = \lVert Au_{\phi_m} \rVert_H$, $X = C\lVert u_{\phi_n} - u_{\phi_m} \rVert_V^2$ and 
$Y = \lVert A(u_{\phi_n} - u_{\phi_m}) \rVert_H^2$ and $Z$ stand for the remaining terms in the right hand side of the equation \ref{main_tau}, that we prove converging to 0. Moreover, $\frac{1}{6}\nu \lVert A(u_{\phi_n} - u_{\phi_m})\rVert_H^2$ is absorbed to the left hand side of the main equation. 

Next, we treat the second nonlinear term as
\begin{align}
&|\langle B(u_{\phi_n}, u_{\phi_n} - u_{\phi_m}), A(u_{\phi_n} - u_{\phi_m}) \rangle|
\nonumber \\& \indeq 
\leq \lVert u_{\phi_n} \rVert_V \lVert u_{\phi_m} -u_{\phi_n} \rVert_V^{1/2} \lVert A(u_{\phi_n} - u_{\phi_m}) \rVert_H^{3/2}
\nonumber \\& \indeq 
\leq \frac{1}{6} \lVert  A(u_{\phi_n} - u_{\phi_m})  \rVert_H^2 + C(M,\tilde{M})\lVert u_{\phi_n} - u_{\phi_m} \rVert_V^2,
\end{align}
where the first term is absorbed to LHS of the equation \ref{main_tau}, whereas for the second term we have 
\begin{equation}
 C(M,\tilde{M})\mathbb{E}[\int_0^{\tau}\lVert u_{\phi_n} - u_{\phi_m} \rVert_V^2] \rightarrow 0,
\end{equation}
as $m,n \rightarrow \infty$ by Theorem 3.5. Hence , the first part of the proof is concluded.

2. The proof is identical with \cite{GZ09}. First by Ito we have 
\begin{align}
d \lVert u_\phi \rVert_V^2 + 2 \nu \lVert Au_\phi \rVert_H^2dt &= \bigg( 2\langle \phi - B(u_\phi), Au_\phi \rangle +  \lVert g_k(u_\phi) \rVert_V^2 \bigg)
\nonumber \\& \indeq 
+ 2 \sum_{k=1}^\infty \langle g_k(u_\phi),Au_\phi \rangle dW_k
\end{align}
We fix $\tau \in \mathcal{T}_n^{M,T}$ and $S > 0$. Integrating from $0$ to $\tau \wedge S$, we get 
\begin{align}
\sup_{r \in  [0,S \wedge \tau]} \lVert u_\phi \rVert_V^2 + \int_0^{S \wedge \tau} 2 \nu \lVert Au_\phi \rVert_H^2 ds &\leq \lVert u_0 \rVert_V^2 + \int_0^{S \wedge \tau} 2 | \langle \phi - B(u_\phi),Au_\phi \rangle | dr 
\nonumber \\& \indeq 
+ \int_0^{S \wedge \tau} \lVert g(u_\phi) \rVert_V^2 dr
\nonumber \\& \indeq 
+ \sup_{r \in [0, S \wedge \tau ]}\bigg| \sum_{k = 1}^\infty 2 \langle g_k(u_\phi),u_\phi \rangle dW_k  \bigg|.
\end{align}
Applying the classical estimate on nonlinear term (see \cite{CF88}), we have 
\begin{align} \label{clas_es} 
\bigg| \langle B(u_\phi),A(u_\phi) \rangle \bigg| &\leq \lVert u_\phi \rVert_V^{3/2} \lVert Au_\phi \rVert_H^{3/2}
\nonumber \\&
\leq C \lVert u_\phi \rVert_V^6 + \frac{\nu}{4}\lVert Au_\phi \rVert_H^2.
\end{align}
Using Equation (\ref{clas_es}) and the Lipschitz assumption on $g$, we get
\begin{align}
\sup_{r \in [0,S \wedge \tau]} \lVert u_\phi \rVert_V^2 &+ \int_0^{\tau \wedge S} \nu \lVert Au_\phi \rVert_H^2 dr
\nonumber \\& \indeq 
\leq \lVert u_0 \rVert_V^2 + C\int_0^{S \wedge \tau} ( \lVert \phi \rVert_H^2 + \lVert u_\phi \rVert_V^6 + \lVert u_\phi \rVert_V^2 +1 )dr
\nonumber \\& \indeq 
+ \sup_{ r \in [0,S \wedge \tau] } \bigg| \int_0^r 2 \sum_{k=1}^\infty \langle g_k(u_\phi),u_\phi \rangle dW_k \bigg|.
\end{align}
This implies then 
\begin{align}
&\mathbb{P}\bigg( \sup_{s \in [0,\tau \wedge S]} \lVert u_\phi \rVert_V^2 + \nu \int_0^{\tau \wedge S} \lVert Au_\phi \rVert_H^2ds > \lVert u_0 \rVert_V^2 + (M-1)^2 \bigg)
\nonumber \\& \indeq
\leq \mathbb{P}\bigg( C \int_0^{\tau \wedge S} (\lVert \phi \rVert_H^2 + \lVert u_\phi \rVert_V^6 + \lVert u_\phi \rVert_V^2 + 1) dr > \frac{(M-1)^2}{2} \bigg)
\nonumber \\& \indeq
\leq \frac{2C}{(M-1)^2}\mathbb{E}\int_0^{\tau \wedge S} (\lVert \phi \rVert_H^2 + \lVert u_\phi \rVert_V^6 + \lVert u_\phi \rVert_V^2 + 1)dr 
\nonumber \\& \indeq
\leq C \mathbb{E}\bigg( \int_0^S \lVert\phi\rVert^2_H + 1 dr \bigg)
\end{align}
Next, using Doob's Inequality for the second term, we get 
\begin{align} 
&\mathbb{P}\bigg( \sup_{r \in [0,\tau \wedge S]} \bigg| \sum_{k=1}^\infty \int_0^r \langle g_k(u_\phi),u_\phi \rangle dW_k \bigg| > \frac{(M-1)^2}{2} \bigg)
\nonumber \\& \indeq
\leq \frac{4}{(M-1)^4} \mathbb{E}\bigg( \int_0^{S \wedge \tau} \lVert u_\phi \rVert_V^2 \sum_{k=1}^\infty \lVert g_k(u_\phi) \rVert_V^2 dr \bigg)
\nonumber \\& \indeq
\leq C
\end{align}
By letting $S \rightarrow 0$ with the integrability condition imposed on function $\phi$, we conclude the proof. 
\end{proof}
\begin{Theorem} Given the assumptions on initial data $u_0$ as in Definition 2.1 and $\phi \in \mathcal{U}^b$,  there exists a global strong solution $(u_\phi,\tau)$ in the sense of Definition 2.1 introduced above.
\end{Theorem}
\begin{proof} Let $w \in H$ be given. Using Theorem \ref{main_prop} with $\{u_{\phi_n}\}$ be the sequence of solutions of (3.19). Due to \ref{cond1} and \ref{cond2}, we apply Lemma \ref{lemma1} with $B_1 = V$ and $B_2 = D(A)$ and the sequence $\{X^n\}= \{u_{\phi_n}\}$. We infer the existence of a subsequence $\{ u_{\phi_n'}\}$ and a strictly positive stopping time $\tau \leq T$ and a process $u(.) = u(. \wedge \tau)$, continuous in $V$ such that 
\begin{equation}  
\label{tag1}
\sup_{t \in [0,\tau]} \lVert u_{\phi_{n'}} - u \rVert^2_V + \nu \int_0^{\tau} \lVert A(u_{\phi_{n'}} - u) \rVert_H^2ds \rightarrow 0, 
\end{equation}
a.s. We also have that the conditions of Lemma \ref{lemma1} (ii) is satisfied for any $p \in (1,\infty)$. Thus, we have for any $p >1$
\begin{equation}
u_\phi(. \wedge \tau) \in L^p(\Omega; C([0,T];V)),
\end{equation}
with 
\begin{equation}
u_\phi\mathbbm{1}_{t \leq \tau} \in L^p(\Omega;L^2([0,T];D(A))).
\end{equation}
By Lemma \ref{lemma1} (ii) \, we infer a collection of measurable sets $\Omega_{n'} \in \mathcal{F}$ with 
\begin{equation}
\label{tag2}
\Omega_{n'} \uparrow \Omega
\end{equation}
such that 
\begin{equation}
\label{tag3}
\sup_{n'} \mathbb{E}\bigg[ \sup_{t \in [0,\tau]} \lVert u_{\phi_{n'}}(t) \mathbbm{1}_{\Omega_{n'}} \rVert_V^2 + \nu \int_0^\tau \lVert Au_{\phi_{n'}}\mathbbm{1}_{\Omega_{n'}} \rVert_H^2 ds \bigg]^{p/2} < \infty.
\end{equation}
Using equations \ref{tag1}, \ref{tag2}, \ref{tag3} and by Lemma \ref{lemma2} we have 
\begin{equation}
\mathbbm{1}_{\Omega_{n'}, t \leq \tau} u_\phi^{n'} \rightharpoonup \mathbbm{1}_{t \leq \tau}u_\phi \textrm{ in } L^p(\Omega;L^2( [0,T];D(A))),
\end{equation}
and
\begin{equation}
\mathbbm{1}_{\Omega_{n'}} u_\phi^{n' \wedge \tau} \rightharpoonup^* u_\phi \textrm{ in } L^p(\Omega;L^\infty( [0,T];V)).
\end{equation}
For the  nonlinear term, we estimate for all $w \in H$ as follows: 
\begin{align}
&\big| \langle B(u_{\phi_{n'}},u_{\phi_{n'}}) - B(u,u),w \rangle \big| 
\\ \nonumber
&= \big| \langle B(u_{\phi_{n'}}-u,u_{\phi_{n'}}) + B(u,u_{\phi_{n'}}) - B(u,u), w \rangle \big|
\\ \nonumber
&= \big| \langle B(u_{\phi_{n'}}-u,u_{\phi_{n'}}) + B(u,u_{\phi_{n'}}-u),w \rangle \big|
\\ \nonumber
&\leq \big| \langle B(u_{\phi_{n'}}-u,u_{\phi_{n'}}),w \rangle  \big| + \big| \langle B(u,u_{\phi_{n'}}-u),w \rangle \big|
\end{align}
Then we have using the classical estimates [CF]
\begin{align}
\big| \langle B(u_{\phi_{n'}}-u,u_{\phi_{n'}}),w \rangle \big| &\leq C \lVert u_{\phi_{n'}} - u \rVert_H^{1/2}\lVert u_{\phi_{n'}} - u \rVert_V^{1/2}\lVert u_{\phi_{n'}} \rVert_V^{1/2}\lVert Au_{\phi_{n'}} \rVert_H^{1/2}\lVert w \rVert_H^{1/2}
\\ \nonumber
\big| \langle B(u,u_{\phi_{n'}} - u),w \rangle \big| &\leq C\lVert u \rVert_H^{1/2}\lVert u \rVert_V^{1/2}\lVert u_{\phi_{n'}} - u \rVert_V^{1/2}\lVert A(u - u_{\phi_{n'}}) \rVert_H^{1/2}\lVert w \rVert_H^{1/2}
\end{align}
Hence the nonlinear terms converge to 0 by \ref{tag1}, we conclude that given any $v \in H$
\begin{equation}
\label{tag4}
\mathbbm{1}_{t \leq \tau}\langle B(u_{n'-1}, u_{n'}),v \rangle \rightarrow \mathbbm{1}_{t \leq \tau}\langle B(u,u), v \rangle,
\end{equation}
as $n' \rightarrow \infty$, for almost every $(\omega,t) \in \Omega \times [0,T]$.
Moreover, by using the uniform bound of Equation \ref{tag1} with $p = 4$, one finds that 
\begin{align}
\label{tag5}
&\sup_{n'}\mathbb{E}\bigg( \mathbbm{1}_{\Omega_{n'}}\int_0^{\tau} | B(u_{\phi_{n'}},u_{\phi_{n'}}) |^2 ds \bigg)
\\ \nonumber \indeq 
&\leq C \sup_{n'}\mathbb{E}\bigg( \mathbbm{1}_{\Omega_{n'}} \int_0^{\tau}  \lVert u_{\phi_{n'}} \rVert_H^2 \lVert u_{\phi_{n'}} \rVert_V^2 ds \bigg)
\\ \nonumber \indeq 
&\leq C \sup_{n'}\mathbb{E}\bigg( \mathbbm{1}_{\Omega_{n'}} \sup_{t \in [0,\tau]}\lVert u_{\phi_{n'}} \rVert_H^3 \lVert Au_{n'} \rVert_H ds \bigg)
\\ \nonumber \indeq 
&\leq C \sup_{n'}\mathbb{E}\bigg( \mathbbm{1}_{\Omega_{n'}} \sup_{t \in [0,\tau]}\lVert u_{\phi_{n'}} \rVert_H^4 + \bigg( \int_0^\tau \lVert Au_{n'} \rVert_H^2 ds \bigg)^2  \bigg)
\\ \nonumber \indeq 
& < \infty
\end{align}
Using \ref{tag4} and \ref{tag5} and using Lemma \ref{lemma2} we have 
\begin{equation}
\mathbbm{1}_{\Omega_{n'},t\leq \tau}B(u_{\phi_{n'}},u_{\phi_{n'}}) \rightharpoonup \mathbbm{1}_{t \leq \tau}B(u,u),
\end{equation}
in $L^2(\Omega;L^2([0,T];H))$.
By Lipschitz condition on $g$ we get 
\begin{align} 
&\sum \lVert g_k(u_{\phi_{n'}}) - g_k(u_\phi) \rVert_V^2 
\\ \nonumber \indeq 
&\leq \lVert u_\phi - u_{\phi_{n'}} \rVert_V^2 \rightarrow 0,
\end{align}
by \ref{tag1}. We have moreover that 
\begin{align}
&\sup_{n'}\mathbb{E} \bigg[ \mathbbm{1}_{\Omega_{n'}} \int_0^{\tau} \lVert g_k(u_{\phi_{n'}}) \rVert_V^2ds \bigg]
\\ \nonumber \indeq 
&\leq C \sup_{n'}\mathbb{E}\bigg[ \mathbbm{1}_{\Omega_{n'}} \int_0^{\tau} 1 + \lVert u_{n'} \rVert_V^2ds \bigg]
\\ \nonumber \indeq 
&< \infty,
\end{align}
which means that 
\begin{equation}
\mathbbm{1}_{\Omega_{n'},t \leq \tau} g(u_{\phi_{n'}}) \rightarrow \mathbbm{1}_{t \leq \tau}g(u_\phi),
\end{equation}
in $L^2(\Omega;L^2( [0,T];l^2(H) ))$.
Next we have 
\begin{align}
&\mathbb{E}[ 2 \int_0^\tau | \langle \phi_n (t, u_{\phi_n}) - \phi_m (t, u_{\phi_m}), A(u_{\phi_n} - u_{\phi_m})  \rangle | dt ]
\nonumber \\& \quad
\leq
\int_0^{\tau} \lVert \phi_n(t,u_{\phi_n}) - \phi_m(t,u_{\phi_m}) \rVert_V \lVert u_{\phi_n} - u_{\phi_m} \rVert_V ds 
\nonumber \\& \quad \quad
+ C \mathbb{E}\int_0^T \lVert u_{\phi_n} - u_{\phi_m} \rVert_V^2
\nonumber \\& \quad
\leq C \mathbb{E}\int_0^\tau \lVert \phi_n(t,u_{\phi_n}) - \phi_n(t,u_{\phi_m})  \rVert_V^2
\nonumber \\& \quad\quad
+ C\mathbb{E}\int_0^\tau \lVert \phi_n(t,u_{\phi_m}) - \phi_m(t,u_{\phi_m}) \rVert_V^2
\nonumber \\& \quad\quad
+ C \mathbb{E}\int_0^T \lVert u_{\phi_n} - u_{\phi_m} \rVert_V^2
\end{align}
which goes to 0 as $n,m \rightarrow \infty$ by Theorem \ref{breckner} as well as by our assumption on $\phi$.

Hence by combining all the estimates and using ${\phi_{n'}}$ being bounded as well we deduce that for any fixed $v \in H$
\begin{align}
\label{eqn100}
\mathbbm{1}_{\Omega_{n'}}\int_0^{t \wedge \tau} \langle Au_{\phi_{n'}},v \rangle ds &\rightharpoonup \int_0^{t \wedge \tau} \langle Au_\phi,v \rangle ds,
\\ \nonumber
\mathbbm{1}_{\Omega_{n'}}\int_0^{t \wedge \tau} \langle B(u_{\phi_{n'}}),v \rangle ds &\rightharpoonup \int_0^{t \wedge \tau} \langle B(u_\phi),v \rangle ds,
\\ \nonumber
\mathbbm{1}_{\Omega_{n'}} \int_0^{t \wedge \tau} \langle \phi_{n'},v \rangle ds &\rightharpoonup \int_0^{t \wedge \tau} \langle \phi,v \rangle ds.
\\ \nonumber
\mathbbm{1}_{\Omega_{n'}} \sum_k \int_0^{t \wedge \tau} g_k(u_{\phi_{n'}},v)dW_k &\rightharpoonup \int_0^{t \wedge \tau} \langle g_k(u_\phi),v\rangle dW_k.
\end{align}
weakly in $L^2(\Omega \times [0,T])$. If $K \subset \Omega \times [0,T]$ is any measurable set, then by \ref{eqn100}, we have 
\begin{align}
&\mathbb{E} \int_0^T \mathbbm{1}_K(\omega, t)\langle u_\phi(t \wedge \tau),v \rangle dt 
\\ \nonumber & \indeq 
= \lim_{n'\rightarrow \infty} \mathbb{E}\int_0^T \langle \mathbbm{1}_{\Omega_{n'}}(\omega), \mathbbm{1}_K(\omega, t)v \rangle dt
\\ \nonumber & \indeq 
= \lim_{n'\rightarrow \infty} \bigg( \mathbb{E}\int_0^T \mathbbm{1}_K(\omega, t)\mathbbm{1}_{\Omega_{n'}}(\omega)\langle u_0,v \rangle dt 
\\ \nonumber & \indeq \indeq 
- \mathbb{E}\int_0^T \mathbbm{1}_K(\omega, t) \mathbbm{1}_{\Omega_{n'}}(\omega)\bigg[ \int_0^{t \wedge \tau} \langle\nu Au_{\phi_{n'}} + B(u_{\phi_{n'}} - \phi_{n'},v \rangle ds \bigg] dt
\\ \nonumber & \indeq \indeq 
+ \mathbb{E}\int_0^T \mathbbm{1}_K(\omega,t)\mathbbm{1}_{\Omega_{n'}}(\omega)\bigg[ \sum_k \int_0^{t \wedge \tau} \langle g_k(u_{\phi_{n'}}),v\rangle  dW_s \bigg]dt \bigg)
\\ \nonumber & \indeq
= \mathbb{E}\int_0^T \mathbbm{1}_K(\omega,t)\bigg[ \langle u_0,v \rangle - \int_0^{t \wedge \tau} \langle \nu Au_\phi + B(u_\phi) -\phi(s,u_{\phi}),v \rangle ds \bigg] dt
\\ \nonumber & \indeq \indeq 
+ \mathbb{E}\int_0^{T} \mathbbm{1}_K(\omega, t) \bigg[ \sum_k \int_0^{t \wedge \tau} \langle g_k(u_\phi),v \rangle dW_k \bigg] dt
\end{align}
Since $v \in H$ and $K$ are arbitrary, we conclude that $u_\phi$ satisfies the regularity conditions. Hence, we have shown the local existence of the solution $u_\phi$. Relaxing the restriction $\lVert u_0 \rVert_V \leq \tilde{M}$, namely extending to the case $u_0 \in L^2(\Omega, V)$ and the global uniqueness follows the same steps of \cite{GZ09} Theorem 4.2. This concludes the proof.
\end{proof}
Next, we borrow the following theorem from \cite{KUZ}.
\begin{Theorem} \label{kuz} \cite{KUZ}
Let $u_\phi,u_0,\phi,g$ be as defined above with the corresponding properties. Then, we have 
\begin{equation}
\mathbb{E}[\sup_{[0,T]}  \log(1 + \lVert u_\phi \rVert_V^2)] \leq C(\phi,g,u_0,T).
\end{equation}
\end{Theorem}
We continue with the following lemma.
\begin{Lemma}\label{V_prob} Given the assumptions on initial data and $\{\phi_n\}_{n\geq 1}$ as in Theorem 3.8, we have that for any deterministic time $T$
\begin{equation}
\mathbb{P} ( \sup_{t \in [0,T]} \lVert u_\phi -  u_{\phi_n} \rVert_V^2 > \delta ) < \epsilon
\end{equation}
for any $n,m \geq N_0$ for some $N_0$, i.e. solutions with different deterministic force $\{u_{\phi_n}\}_{n \geq 1}$, converge in probability to $u_\phi$ as $n,m \rightarrow \infty$.
\end{Lemma}
\begin{proof}
By assumption, we have $u_0 \in L^2(\Omega,V)$. Hence, by Chebyshev theorem we have, 
\begin{equation}
\mathbb{P}(\lVert u_0 \rVert_V^2 >s)\rightarrow 0
\end{equation}
as $s \rightarrow \infty$. Denoting $\Omega_s = \{ \lVert u_0 \rVert_V^2 \leq s \}$, we have $\Omega_s \rightarrow \Omega$. Hence, we choose $s$ such that $\mathbb{P}(\Omega_s) > 1 - \frac{\epsilon}{2}$. Moreover, we know by Theorem \ref{main_prop} and Lemma \ref{lemma1} that there exists a sequence of stopping times $\{\tilde{\tau}_{n_l}^M\}_{n_l \geq 1}$ with the corresponding subsequence $\{u_{\phi_{n_l}}\}$ converging monotone decreasing to $\tau^M$ by $\ref{5.12}$. We also know by Theorem 3.7 that $\tau^M\rightarrow \infty$ a.s. as $M\rightarrow \infty$, where $M$ is the constant defined as in Lemma 4.1., since the solution is global in the sense of Definition \ref{global_def}.

Hence, denoting $\{ \tau_{n_l}^M = \tilde{\tau}_{n_l}^M \wedge T \}_{n_l \geq 1}$, there exists $M_0$ such that $\mathbb{P}(\tau^{M_0} < T ) \leq \frac{\epsilon}{4}$ and by Lebesgue dominated convergence theorem, we have
\begin{equation}
\lim_{n_l \rightarrow \infty}\mathbb{E}\left[ \mathbbm{1}_{\Omega_s} \sup_{t \in [0,\tau^{M_0}]} \lVert u_\phi -  u_{\phi_{n_l}} \rVert_V^2\right] = 0.
\end{equation}
This implies convergence in probability. Thus,
\begin{equation}
\lim_{n_l \rightarrow \infty} \mathbb{P} ( \mathbbm{1}_{\Omega_s} \sup_{t \in [0,\tau^{M_0}]} \lVert u_\phi -  u_{\phi_{n_l}} \rVert_V^2 > \delta ) = 0,
\end{equation}
for any $\delta > 0$. Hence, we have
\begin{align}
&\mathbb{P}\big(\mathbbm{1}_{\Omega_s}\sup_{t \in [0,T]} \lVert u_\phi -  u_{\phi_{n_l}} \rVert_V^2 \geq \delta\big) = \mathbb{P}\big(\{\sup_{t \in [0,T]} \lVert u_\phi -  u_{\phi_{n_l}} \rVert_V^2 \geq \delta\}  \cap  \{\tau^{M_0} < T\} \cap  \{\omega \in \Omega_s\} \big) 
\nonumber \\& \indeq
+ \mathbb{P}(\{\sup_{t \in [0,T]} \lVert u_\phi -  u_{\phi_{n_l}} \rVert_V^2 \geq \delta \}  \cap  \{\tau^{M_0} = T\}  \cap  \{\omega \in \Omega_s\} )
\nonumber \\&
\leq \mathbb{P}( \tau^{M_0} < T ) + \mathbb{P}(\mathbbm{1}_{\Omega_s}\sup_{t \in [0,\tau^{M_0}]} \lVert u_\phi -  u_{\phi_{n_l}} \rVert_V^2 > \delta )
\end{align}
Then, we get 
\begin{align}
\mathbb{P}\left( \sup_{t \in [0,T]} \lVert u_\phi - u_{\phi_{n_l}} \rVert_V^2 \geq \delta \right) &\leq \mathbb{P} \left( \tau^{M_0} < T \right) + \mathbb{P} \left( \mathbbm{1}_{\Omega_s}\sup_{t \in [0,\tau^{M_0}]} \lVert u_\phi - u_{\phi_{n_l}} \rVert_V^2 > \delta \right) + \mathbb{P}\left( \Omega_s^c \right)
\nonumber \\&\leq
\frac{\epsilon}{4} + \frac{\epsilon}{4} + \frac{\epsilon}{2} = \epsilon.
\end{align}
for $n_l$ large enough. Then by taking any subsequence $u_{\phi_{m_l}}$ and by Theorem \ref{main_prop} and Lemma \ref{lemma1} repeating the same arguments above, we get that every subsequence $\{u_{\phi_{m_l}}\}$ has a further subsequence that converges in probability to $u_\phi$, which implies that the whole sequence $\{u_{\phi_n}\}$ converges in probability to $u_\phi$, which concludes the proof.
\end{proof}
\newline
Now we are ready to prove Theorem \ref{thm1}. 
\begin{proof}
\begin{align}
\label{main_proof}
&\big|\mathbb{E}\big[\sup_{[0,T]} \varphi (\mathcal{L}(t,u_{\phi_n}, \phi_n)) - \sup_{[0,T]}\varphi (\mathcal{L}(t,u_{\phi}, \phi))\big] \big| 
\nonumber \\& \indeq 
\leq \mathbb{E}\big[\sup_{[0,T]}  \big|\varphi (\mathcal{L}(t,u_{\phi_n}, \phi_n)) - \varphi (\mathcal{L}(t,u_{\phi}, \phi))\big| \big]
\nonumber \\& \indeq 
\leq \mathbb{E}\big[\sup_{[0,T]}  \big|\varphi (\mathcal{L}(t,u_{\phi_n}, \phi_n)) - \varphi (\mathcal{L}(t,u_{\phi_n}, \phi))\big| \big] 
\nonumber \\& \indeq \indeq 
+ \mathbb{E}\big[\sup_{[0,T]}  \big|\varphi (\mathcal{L}(t,u_{\phi_n}, \phi)) - \varphi (\mathcal{L}(t,u_{\phi}, \phi))\big| \big]
\nonumber \\& \indeq
\leq \mathbb{E}\big[\sup_{[0,T]}\varphi (\big|\mathcal{L}(t,u_{\phi_n}, \phi_n)) - \mathcal{L}(t,u_{\phi_n}, \phi) \big| \big] 
\nonumber \\& \indeq \indeq 
+ \mathbb{E}\big[\sup_{[0,T]}\varphi (\big|\mathcal{L}(t,u_\phi, \phi_n)) - \mathcal{L}(t,u_\phi, \phi) \big| \big] 
\nonumber \\& \indeq
\leq \mathbb{E} \big[ \sup_{[0,T]} \varphi (  C\lVert \phi - \phi_n \rVert_H^2 + C\lVert u_{\phi} - u_{\phi_n} \rVert_V^2 )  \big]
\nonumber \\& \indeq 
\leq \mathbb{E} \big[  \sup_{[0,T]} \varphi( C \lVert \phi - \phi_n \rVert_H^2)  \big] + \mathbb{E} \big[ \sup_{[0,T]} \varphi ( C\lVert u_{\phi} - u_{\phi_n} \rVert_V^2)  \big],
\end{align}
where we appeal to  Lemma \ref{con_lem} in the third inequality and the Lipschitz assumption on $\mathcal{L}(t,u_{\phi}, \phi)$ in the fourth inequality. 
We have by boundedness assumption on $\{ \phi_n\}$ the followings
\begin{align} 
\sup_{[0,T]} \varphi( C \lVert \phi - \phi_n \rVert_H^2) \rightarrow 0, \textrm{ as  } n &\rightarrow \infty
\nonumber \\
\mathbb{E}[ \sup_{[0,T]}  \lVert \phi - \phi_n \rVert_H^2 ] \leq M,
\end{align}
Hence, we have by uniform integrability and convergence in probability
\begin{equation}
\mathbb{E}[ \sup_{[0,T]} \varphi( \lVert \phi - \phi_n \rVert_H^2 ) ] \rightarrow 0, \textrm{ as } n \rightarrow \infty.
\end{equation}
For the first term in the last line of the Equation \ref{main_proof}, we have 
\begin{equation} 
\mathbb{E} \big[ \sup_{[0,T]} \varphi ( C\lVert u_{\phi} - u_{\phi_n} \rVert_V^2) \big] \leq \mathbb{E} \big[ \sup_{[0,T]} \varphi(C\lVert u_{\phi}\rVert_V^2 ) \big] + \mathbb{E} \big[ \sup_{[0,T]} \varphi(C\lVert u_{\phi_n}\rVert_V^2 ) \big]
\end{equation}
By noting that $\log(1+x) \leq x$ for $x > 0$, we have by Lemma \ref{V_prob} that 
\begin{equation}
\mathbb{P}\left(\sup_{[0,T]} (  \log(1 + \lVert u_{\phi} - u_{\phi_n} \rVert_V^2)^{1-\epsilon}\right) \rightarrow 0
\end{equation}
for $0 < \epsilon < 1$ in probability as $n \rightarrow \infty$. Moreover, using Theorem \ref{kuz}, we have that 
\begin{equation}
\mathbb{E}\left[ \left(\sup_{[0,T]} (\log(1 + \lVert u - u^n \rVert_V^2)\right)^{(1-\epsilon){\frac{1}{1-\epsilon}}} \right] \leq M(u_0,f,g,T).
\end{equation}
We note here that $g(x) = x^{\frac{1}{1-\epsilon}}$ is a convex function with $\lim_{x\rightarrow \infty} \frac{x^{\frac{1}{1-\epsilon}}}{x} = \infty$. Using de La-Vallee-Poussin criteria for uniform integrability (see e.g. \cite{D13}) we get that
\begin{equation}
\left \{ \sup_{[0,T]} \left( \log(1 + \lVert u_{\phi} - u_{\phi_n} \rVert_V^2 \right)^{1-\epsilon} \right \}_{n \geq 1}
\end{equation} is uniformly integrable. Using uniform integrability and by Lemma \ref{V_prob} convergence in probability imply $L^{1}$-convergence \cite{D13}. Thus, using that $x^{1-\epsilon}$ for $0 < \epsilon < 1$ being increasing and continuous, we get that 
\begin{equation}
\mathbb{E}\left[ \sup_{[0,T]} \left(\log(1 + \lVert u_{\phi} - u_{\phi_n} \rVert_V^2 \right)^{1-\epsilon} \right] \rightarrow 0
\end{equation}
as $n \rightarrow \infty$. Hence, Theorem \ref{thm1} is proven. 
\end{proof}
By generalized Weierstrass theorem (see \cite{Z85}), we conclude the following corollary.
\begin{Corollary} \label{cor2}
In addition to the assumptions specified above, if the admissible controls are in a compact subset of $\mathcal{U}$, then there exists an optimal feedback control $\phi^*$ satisfying 
\begin{equation}
J(\phi^*) = \min_{\phi \in \mathcal{U}} J(\phi)
\end{equation}
\end{Corollary}

\section{Appendix}
\begin{Lemma} \cite{GZ09}
\label{lemma1}
Let $(\Omega,\mathcal{F},(\mathcal{F}_t)_{t\geq 0},\mathbb{P})$ be a fixed filtered probability space. Suppose that $B_1$ and $B_2$ are Banach spaces with $B_2 \subset B_1$ with continuous embedding. We denote the associated norms by $|\cdot|_i$. Define 
\begin{equation}
\mathcal{E}(T) := C([0,T];B_1) \cap L^2([0,T];B_2)
\end{equation}
with the norm 
\begin{equation}
|Y|_{\mathcal{E}(T)} = \bigg( \sup_{t \in [0,T]}|Y(t)|_1^2 + \int_0^T |Y(t)|_2^2dt \bigg)^{1/2}.
\end{equation}
Let $X_n$ be a sequence of $B_2$-valued stochastic process such that for every $T>0$ we have $X_n \in \mathcal{E}(T)$, a.s. For $M > 1$, $T > 0$ define the collection of stopping times 
\begin{equation}
\mathcal{T}_n^{M,T} := \{ \tau \leq T: |X_n|_{\mathcal{E}(\tau)} \leq M + |X_n(0)|_1 \},
\end{equation}
and let $\mathcal{T}_{n,m}^{M,T} : \mathcal{T}_n^{M,T} \cap \mathcal{T}_m^{M,T}$.
\begin{enumerate}
\item Suppose that for $M > 1$ and $T$, we have 
\begin{equation}
\label{5.3}
\lim_{n \rightarrow \infty}\sup_{m \geq n}\sup_{\tau in \mathcal{T}_n^{M,T}}\mathbb{E}|X_n - X_m|_{\mathcal{E}(\tau)} = 0
\end{equation}
and 
\begin{equation}
\label{5.4}
\lim_{S \rightarrow 0}\sup_n\sup_{\tau \in \mathcal{T}_n^{M,T}}\mathbb{P}[ |X_n|_{\mathcal{E}(\tau \wedge S)} > |X_n(0)|_1 + M-1 ] = 0.
\end{equation}
Then, there exists a stopping time $\tau$ with:
\begin{equation}
\label{5.5}
\mathbb{P}(0 < \tau \leq T) = 1,
\end{equation}
and a process $X(\cdot) = X(\cdot \wedge \tau \in \mathcal{E}(\tau))$, such that 
\begin{equation}
\label{5.6}
|X_{n_l} - X|_{\mathcal{E}(\tau)} \rightarrow 0, \textrm{ a.s. }
\end{equation}
for some subsequence $n_l \uparrow \infty$. Moreover
\begin{equation}
\label{5.7}
|X|_{\mathcal{E}(\tau)} \leq M + \sup_n |X_n(0)|_1, \textrm{ a.s. }
\end{equation}
\item If, in addition to the conditions imposed above, we also have 
\begin{equation}
\label{5.8}
\sup_n \mathbb{E}|X_n(0)|_1^p < \infty,
\end{equation}
for some $1 \leq p < \infty$, then there exists a sequence of sets $\Omega_l \uparrow \Omega$ such that 
\begin{equation}
\label{5.9}
\sup_l \mathbb{E}I_{\Omega_l}| X_{n_l} |_{\mathcal{E}(\tau)}^p < \infty
\end{equation}
and 
\begin{equation}
\label{5.10}
\mathbb{E}|X|_{\mathcal{E}(\tau)}^p \leq C_q(M^p + \sup_n \mathbb{E}|X_n(0)|_1^p).
\end{equation}
\end{enumerate}
\end{Lemma}
\begin{proof} To find the convergent subsequence, we proceed by induction on $l$ and start with $l=0$ and $n_0 =1$. We have by \ref{5.3}
\begin{equation}
\label{5.11}
  \sup_{\tau  \in \mathcal{T}_{n_{l+1},n_l}^{M,T}} \mathbb{E}| X_{n_l} - X_{n_{l+1}} |_{\mathcal{E}(\tau)} \leq 2^{-2l}.
\end{equation}
Next to find $\tau$ in \ref{5.5} and \ref{5.6}, we define 
\begin{equation}
\label{5.12}
\tau_l := \inf_{t > 0}\{ | X_{n_l} |_{\mathcal{E}(t)} > | X_{n_l}(0) |_1 + (M - 1 + 2^{-l}) \} \wedge T, 
\end{equation}
and let 
\begin{equation}
\label{5.15}
\Omega_N = \bigcap_{j = N}^\infty \bigg\{  | X_{n_j} - X_{n_{j+1}} |_{\mathcal{E}(\tau_j \wedge \tau_{j+1})} < 2^{-(l-2)} \bigg\}
\end{equation}
Using $\tau_l \wedge \tau_{l+1} \in \mathcal{T}_{n_{l+1},n_l}^{M,T}$, we have
\begin{equation}
\mathbb{P}\big( |X_{n_l} - X_{n_{l+1}} |_{\mathcal{E}(\tau_l \wedge \tau_{l+1})} \geq 2^{-(l+2)} \big) \leq 2^{l+2}\mathbb{E}| X_{n_l} - X_{n_{l+1}} |_{\mathcal{E}(\tau_l \wedge \tau_{l+1})} \leq 2^{-(l-2)},
\end{equation}
Hence, by Borel-Cantelli lemma, we conclude that
\begin{equation}
\label{5.16}
\mathbb{P}\bigg( \bigcap_{N=1}^\infty \bigcup_{j=N}^\infty \big\{ | X_{n_j} - X_{n_{j+1}} |\mathcal{E}(\tau_j \wedge \tau_{j+1}) \geq 2^{-(j+2)} \big\}  \bigg) = 0
\end{equation}
and hence $\tilde{\Omega} := \cup_N \Omega_N$ is a set of full measure. Next, we note that 
\begin{equation}
\label{5.17}
\tau_{l+1}(\omega) \leq \tau_l(\omega),
\end{equation}
for every $l \geq N, \omega \in \Omega_N$, since given $N$ and $l \geq N$, take the set $\{ \tau_{l+1} > \tau_l \} \cap \Omega_N$. On this set, we have $\tau_l < T$. By continuity of $|X_{n_l}|_{\mathcal{E}(t)}$ in t, that imples
\begin{equation}
|X_{n_l}|_{\mathcal{E}(\tau_l)} = |X_{n_l}(0)|_1 + (M - 1 + 2^{-l}).
\end{equation}
Moreover, we have on $\Omega_N$
\begin{align}
|X_{n_l}|_{\mathcal{E}(\tau_l \wedge \tau_{l+1})} - |X_{n_{l+1}}|_{\mathcal{E}(\tau_l \wedge \tau_{l+1})} &< 2^{-(l+2)}
\\ \nonumber
|X_{n_{l+1}}(0)| - |X_{n_l}(0)| < 2^{-(l+2)}.
\end{align}
Hence, we get that 
\begin{align}
\label{5.18}
|X_{n_{l+1}}|_{\mathcal{E}(\tau_l \wedge \tau_{l+1})} &> |X_{n_l}|_{\mathcal{E}(\tau_l \wedge \tau_{l+1})} - 2^{-(l+2)}
\\ \nonumber
&= |X_{n_l}|_{\mathcal{E}(\tau_l)} - 2^{-(l+2)}
\\ \nonumber
&= |X_{n_l}(0)|_1 + (M - 1 + 2^{-l}) - 2^{-(l+2)}
\\ \nonumber
&> |X_{n_{l+1}}(0)|_1 + (M - 1 + 2^{-l}) - 2\cdot2^{-(l+2)}
\\ \nonumber
&= |X_{n_{l+1}(0)}|_1 + (M - 1 + 2^{-l+1}),
\end{align}
over $\{ \tau_{l+1} > \tau_l \} \cap \Omega_N$. Moreover on $\Omega_N$, we have 
\begin{align}
\label{5.19}
|X_{n_{l+1}}|_{\mathcal{E}(\tau_l \wedge \tau_{l+1})} &\leq |X_{n_{l+1}}|_{\mathcal{E}(\tau_{l+1})}
\\ \nonumber \qquad
& \leq |X_{n_{l+1}}(0)| + (M - 1 + 2^{-(l+1)}).
\end{align}
Hence, we have by \ref{5.18} and \ref{5.19} that $\{ \tau_{l+1} > \tau_l \cap \Omega_N \}$ is empty. Hence, by \ref{5.16} and \ref{5.17}, we have 
\begin{equation}
    \tau = \lim_l \tau_l,\textrm{ a.s.}
\end{equation}
Next, by fixing $\epsilon > 0$ with $T > \epsilon > 0$. We have 
\begin{align}
\{ \tau_l < \epsilon \} &\subset \{ |X_{n_l}|_{\mathcal{E}(\tau_l \wedge \epsilon)} = |X_{n_l}(0)|_1 + ( M - 1 + 2^{-l}) \}
\\ \nonumber & \qquad
\subset \{ |X_{n_l}|_{\mathcal{E}( \tau_l \wedge \epsilon )} > |X_{n_l}(0)|_1 + (M-1) \}.
\end{align}
Since
\begin{align} 
\mathbb{P}(\tau < \epsilon) &= \mathbb{P}\bigg( \bigcap_{l=1}^{\infty} \bigcup_{k=l}^\infty \{ \tau_k < \epsilon \} \bigg)
\\ \nonumber  
&\leq \lim\sup_l \mathbb{P}(\tau_l < \epsilon)
\\ \nonumber
&\leq \sup_l \mathbb{P}( |X_{n_l}|_{\mathcal{E}(\tau_l \wedge \epsilon)} > | X_{n_l}(0) |_1 + (M-1) ),
\end{align} by \ref{5.4}, we have 
\begin{equation}
\mathbb{P}(\tau = 0) = \mathbb{P}(\cap_{\epsilon > 0} \{ \tau < \epsilon \}) = \lim_{\epsilon \downarrow 0}\mathbb{P}(\tau < \epsilon) = 0.
\end{equation}
Hence, $\tau \leq T$. Next, we show that $X_{n_l}$ is Cauchy in $\mathcal{E}(\tau)$ a.s. By \ref{5.15}, for every $\omega \in \tilde{\Omega}$, we can choose $N = N(\omega)$ so that $\omega \in \Omega_N$ and $\tau(w) \leq \tau_{l+1}(w) \leq \tau_l(w)$ whenever $l \geq N$. Hence 
\begin{equation}
| X_{n_l}(\omega) - X |_{\mathcal{E}(\tau)} = 0,\textrm{ a.s. }
\end{equation}
Hence , the first part is proven.  Next, we take $\Omega_l$ as in 5.15. Hence by 
\begin{align}
\label{5.30}
\mathbbm{1}_{\Omega_l}|X_{n_l}|_{\mathcal{E}(\tau)} &\leq 2^{-(l+2)} + \mathbbm{1}_{\Omega_l}| X_{n_{l+1}} |_{\mathcal{E}(\tau)}
\\ \nonumber
&\leq | X_{n_{l+1}}(0) |_1 + M
\\ \nonumber
&\leq \sup_n|X_n(0)|_1 + M,
\end{align}
which implies \ref{5.9}. Moreover, by \ref{5.30} we have 
\begin{equation}
\label{5.31}
\mathbb{E}\big( \mathbbm{1}_{\Omega_l}| X_{n_l} |_{\mathcal{E}(\tau)}^p \big) \leq C_p(M^p + \mathbb{E}|X_{n_l}(0)|_1^p)
\end{equation}
Combining \ref{5.10} and \ref{5.31}, we conclude the bound in \ref{5.11}. By Fatou's lemma, we conclude the result in \ref{5.10}. Hence, the proof is concluded. 
\end{proof}
\begin{Lemma}\cite{GZ09} \label{lemma2} Suppose that $X$ is a seperable Banach space and let $D \subset X$ be a dense subset. Let $X^*$ be the dual of $X$ and denote the dual pairing between $X$ and $X^*$ by $\langle \cdot,\cdot \rangle$. Assume that $(E,\mathcal{E},\mu)$ is a finite measure space and that $p \in (1,\infty)$. Assume that $u,u_n \in L^p(E,X^*)$ with $\{u^n\}$ uniformly bounded in $L^p(E,X^*)$ and 
\begin{equation}
\langle u^n,y \rangle \rightarrow \langle u,y \rangle \textrm{  }\mu-a.e.
\end{equation}
for all $y \in D$. Then 
\begin{equation}
u^n \rightharpoonup^* u
\end{equation}
in $L^p(E,X^*)$.
\end{Lemma}

\end{document}